\documentclass{amsart}
\usepackage{amssymb,latexsym}

\newtheorem{theorem}{Theorem}
\newtheorem{proposition}[theorem]{Proposition}
\newtheorem{corollary}[theorem]{Corollary}

\newtheorem{remark}[theorem]{Remark}
\newtheorem{lemma}[theorem]{Lemma}

\def\proof{\smallskip\noindent {\it Proof --- \ }}
\def\proofof#1{\smallskip\noindent {\it Proof of #1 --- \ }}
\def\endproof{\hfill$\square$\medskip}

\def\ind{\mathrm{ind}}

\def\<{\left<}
\def\>{\right>}

\def\C{\mathbb{C}}
\def\CC{\mathrm{C}}
\def\cC{\mathcal{C}}

\def\H{\mathrm{H}}
\def\I{\mathcal{I}}
\def\tI{\widetilde{\mathcal{I}}}

\def\Ker{\mathrm{Ker}}
\def\Ess{\mathrm{Ess}}
\def\R{\mathbb{R}}
\def\S{\mathcal{S}}
\def\Sym{\mathrm{Sym}}
\def\SL{\mathrm{SL}}

\def\Z{\mathbb{Z}}

\def\act{\mathrm{act}}
\def\g{\mathfrak{g}}
\def\h{\mathfrak{h}}

\def\inv{\mathrm{inv}}

\def\pa{\partial}
\def\ttheta{\widetilde{\theta}}

\title[Algebras of Curvature Forms on Homogeneous Manifolds]
{Algebras of Curvature Forms\\ on Homogeneous Manifolds}

\author{Alexander Postnikov \and Boris Shapiro \and Mikhail Shapiro}

\address{Dept.\ Math., Massachusetts Institute of Technology, 
         Cambridge, MA 02139, U.S.A.}
\email{apost@math.mit.edu}
\address{Dept.\ Math., University of Stockholm, Stockholm, S-10691, Sweden}
\email{shapiro@matematik.su.se}
\address{Dept.\ Math., Royal Institute of Technology, Stockholm, 
        S-10044, Sweden}
\email{mshapiro@math.kth.se}

\date{January~17, 1999}

\begin{document}

\maketitle

\centerline{{\footnotesize To Dmitry Borisovich Fuchs with love from former
students of Moscow `Jewish' University}}

\begin{abstract}
Let~$\CC(X)$ be the algebra generated by the curvature two-forms of standard
holomorphic hermitian line bundles over the complex homogeneous manifold
$X=G/B$.
The cohomology ring of~$X$  is a quotient of~$\CC(X)$.  We calculate
the Hilbert polynomial of this algebra.  In particular, we show that the
dimension of~$\CC(X)$ is equal to the number of independent subsets of roots in
the corresponding root system.
We also construct a more general algebra associated with a point
on a Grassmannian.  We calculate its Hilbert polynomial and present the
algebra in terms of generators and relations.
\end{abstract}

\section{Homogeneous Manifolds}

In this section we remind the reader the basic notions and notation related to
homogeneous manifolds~$G/B$ and root systems, as well as fix our terminology.
\medskip

Let~$G$ be a connected complex semisimple Lie group and~$B$ its Borel
subgroup.  The quotient space~$X=G/B$ is then a compact homogeneous complex
manifold.  We choose a maximal compact subgroup~$K$ of~$G$ and denote by
$T=K\cap B$ its maximal torus.  The group~$K$ acts transitively on~$X$.
Thus~$X$ can be identified with the quotient space~$K/T$.
\medskip

By~$\g$ we denote the Lie algebra of~$G$ and by~$\h\subset\g$ its Cartan
subalgebra.   Also denote by~$\g_\R\subset\g$ the real form of~$\g$ such
that $i\,\g_\R$ is the Lie algebra of $K$.  Analogously,
$\h_\R=\h\cap\g_\R$ and
$i\,\h_\R$ is the Lie algebra of the maximal torus $T$.
The {\it root system} associated with~$\g$ is the set~$\Delta$ of nonzero
vectors (roots)~$\alpha\in\h^*$ for which the root spaces
$$
\g_\alpha = \{x\in\g \mid [h,x]=\alpha(h) x\textrm{ for all }h\in\h\}
$$
are nontrivial.  Then $\g$ decomposes into the direct sum of subspaces
$$
\g = \h\oplus\sum_{\alpha\in\Delta} \g_\alpha.
$$
For $\alpha\in\Delta$, the spaces~$\g_\alpha$ and $\h_\alpha=
[\g_\alpha,\g_{-\alpha}]$ are one-dimensional and there exists a unique element
$h_\alpha\in\h_\alpha$ such that $\alpha(h_\alpha)=2$.  The elements
$h_\alpha\in\h$ are called {\it coroots.}  Actually, $\alpha\in\h_\R^*$ and
$h_\alpha\in
\h_\R$, for $\alpha\in\Delta$.  Let us choose generators
$e_\alpha\in\g_\R$ of the root spaces $\g_\alpha$ such that
$[e_\alpha,e_{-\alpha}]=h_\alpha$ for any root $\alpha$.
Then $[h_\alpha,e_\alpha]=2e_\alpha$ and
$[h_\alpha,e_{-\alpha}]=-2e_{-\alpha}$.

The root system~$\Delta$ is subdivided into a disjoint union of sets of
positive
roots~$\Delta_+$ and negative roots~$\Delta_-=-\Delta_+$ such that the direct
sum $\mathfrak{b}=\h\oplus\sum_{\alpha\in\Delta_+}\g_\alpha$ is the Lie algebra
of Borel subgroup~$B$.
The Weyl group~$W$ is the group generated by the
reflections~$s_\alpha:\h^*\to\h^*$,
$\alpha\in\Delta_+$, given by
$$
s_\alpha: \lambda \mapsto \lambda -  \lambda(h_\alpha)\, \alpha.
$$

The lattice $\hat{T}=\{\lambda\in\h^*\mid \lambda(h_\alpha)\in\Z
\textrm{ for all } \alpha\in\Delta\}$ is called the {\it weight lattice.}
Every weight $\lambda\in\hat{T}$ determines an irreducible unitary
representations~$\pi_\lambda: T\to\C^\times$ of the maximal torus~$T$ given by
$\pi_\lambda(\exp(x))=e^{\lambda(x)}$, for $x\in i\,\h_\R$, and every
irreducible
unitary representation of~$T$ is of this form.

For a weight~$\lambda\in\hat{T}$, the
homomorphism~$\pi_\lambda:T\to\C^\times$
extends uniquely to a holomorphic homomorphism~$\bar\pi_\lambda:
B\to\C^\times$.  Thus any~$\lambda\in\hat{T}$ determines a holomorphic line
bundle~$L_\lambda = G\times_B \C = K\times_T\C$ over $X=G/B=K/T$.
The line bundle~$L_\lambda$ has a canonical $K$-invariant hermitian metric.

The classical Borel's theorem~\cite{borel} describes the cohomology
ring~$\H^*(X,\C)$ of the homogeneous manifold~$X$ in terms of generators
and relations:
$$
\H^*(X,\C)\cong \Sym(\h^*)/ I_W,
$$
where~$I_W$ is the ideal in the symmetric algebra~$\Sym(\h^*)$ generated by
the~$W$-invariant elements without constant term.  The natural projection
from~$\Sym(\h^*)$ to~$\H^*(X,\C)$ is the homomorphism that sends a weight
$\lambda\in\hat{T}$ to the first Chern class~$c_1(L_\lambda)$ of the
line bundle~$L_\lambda$.

\medskip
The purpose of this article is to extend the cohomology ring~$\H^*(X,\C)$
to the level of differential forms on~$X$.
It is possible to exhibit differential two-forms that represent the Chern
classes~$c_1(L_\lambda)$ in the de~Rham cohomology of the homogeneous
manifold~$X$.
Recall that for a holomorphic hermitian line bundle~$L:E\to X$ there is
a canonically associated connection on~$E$.  Denote by $\Theta(L)$  the {\it
curvature form} of this connection, which is a differential two-form on~$X$.
Then the form~$i\,\Theta(L)/2\pi$ represents~$c_1(L)$.

In order to construct the curvature forms~$\Theta(L_\lambda)$ explicitly, we
define the elements~$e^\alpha\in\g_\R^*$, $\alpha\in\Delta$, by
$e^{\alpha}(\h)=0$ and $e^\alpha(e_\beta)=\delta_{\alpha\beta}$, for any
$\beta\in\Delta$.  (Here~$\delta_{\alpha\beta}$ is Kroneker's delta.)  The
space of left $K$-invariant differential one-forms on~$K$
can be identified with the dual to its
Lie algebra, t.e., with~$i\,\g_\R^*$.  Thus the elements $i\,e^\alpha$ can be
regarded as one-forms on~$K$.  The differential two-form on~$K$ given by
$$
\phi_\alpha=e^\alpha\wedge e^{-\alpha},\quad \alpha\in\Delta,
$$
is invariant with respect to the right translation action of the torus~$T$.
Thus~$\phi_\alpha$ produces a two-form on the manifold~$X$, for which we
will use the same notation~$\phi_\alpha$.
It is clear from the definition that $\phi_{-\alpha}=-\phi_\alpha$.

The following statement is implicit in~\cite{GS}.

\begin{proposition}  For $\lambda\in\hat{T}$, the curvature form
of the holomorphic hermitian line bundle~$L_\lambda$ is given by
$$
\Theta(L_\lambda) = \sum_{\alpha\in\Delta_+} \lambda(h_\alpha)\,\phi_\alpha.
$$
\label{pr:curvature}
\end{proposition}

Let~$\Phi$ be the algebra generated by the two-forms~$\phi_\alpha$,
$\alpha\in\Delta_+$.
The relations in~$\Phi$ are relatively simple:
$$
\phi_\alpha \phi_\beta = \phi_\beta \phi_\alpha,\\
\qquad\quad(\phi_\alpha)^2=0.
$$
Thus $\Phi$ is a $2^N$-dimensional algebra, where~$N=|\Delta_+|$.
The main object in this paper is the subalgebra of~$\Phi$ generated
by the curvature forms~$\Theta(L_\lambda)$.

\section{Main Results}

Denote by~$\CC(X)$ the subalgebra in the algebra of differential forms on~$X$
that is generated by the curvature forms~$\Theta(L_\lambda)$ of line bundles.
Obviously, $\CC(X)$ has the structure of a graded ring:
$\CC(X)=\CC^0(X)\oplus\CC^1(X)\oplus \CC^2(X)\oplus\cdots$,
where $\CC^k(X)$ is the subspace of~$2k$-forms
in~$\CC(X)$.  In order to formulate our main results about~$\CC(X)$ we need
some extra notation from the matroid theory.

Let~$V$ be a collection of vectors~$v_1,v_2,\dots,v_N$ in a vector space~$E$,
say, over~$\C$.  A subset of vectors in $V$ is called {\it independent} if they
are linearly independent in~$E$.  By convention, the empty subset is 
independent.  Let~$\ind(V)$ be number of all independent subsets in~$V$.

A {\it cycle\/} is a minimal by inclusion not independent subset.
For a cycle $C=\{v_{i_1},\dots,v_{i_l}\}$, there is a unique, up to a factor,
linear dependence $a_1 v_{i_1}+\dots +a_l v_{i_l}=0$
with non-zero $a_i$'s.  Let us fix a linear order
$v_1<v_2<\cdots<v_N$ of all elements of~$V$.
For an independent subset $S$ in $V$, a vector~$v\in V\setminus S$,
is called {\it externally active\/} if the set $S\cup \{v\}$ contains a cycle
$C$ and $v$ is the minimal element of~$C$.  Let $\act(S)$ be the number of
externally active vectors with respect to $S$.

\begin{theorem}
\label{th:main}
The dimension of the algebra $\CC(X)$ is equal to the number
$\ind(\Delta_+)$ of
independent subsets in the set of positive roots~$\Delta_+$.  Moreover, the
dimension of the $k$-th component $\CC^k(X)$ is equal to the number of
independent subsets~$S\subset\Delta_+$ such that $k=N-|S|-\act(S)$,
where~$N=|\Delta_+|$.
\end{theorem}

We remark here that, although the number $\act(S)$ of externally active vectors
depends upon a particular order of elements in~$V$, the total number of
subsets~$S$ with fixed~$|S|+\act(S)$ does not depend upon a choice of ordering.

\bigskip
We will actually prove a more general result about an arbitrary collection of
vectors~$V$.  
Let $V$ and $E$ be as above.  We will assume that the elements 
$v_1,\dots,v_N$ of $V$ span the $n$-dimensional space $E$.  Thus $N\geq n$.
Let $F=\C^N$ be the linear space
with a distinguished basis $\phi^1,\dots,\phi^N$.
Then $V$ defines the projection map $p:F\to E$ that sends the 
$i$th basis element $\phi^i$ to $v_i$.  
The dual map $p^*:E^*\to F^*$ defines an $n$-dimensional
plane $P=\mathrm{Im}(p^*)$ in $F^*$.  In other words, the collection
of vectors $V$ can be identified with an element $P$ of 
the Grassmannian $G(n,N)$ of $n$-dimensional planes in $\C^N$.

Let $\phi_1,\dots,\phi_N$ be the basis in $F^*$ dual to the chosen
basis in $F$.  Denote by $\Phi_N$ the quotient of the 
symmetric algebra $\mathrm{Sym}(F^*)$ modulo the relations $(\phi_i)^2=0$,
$i=1,\dots,N$.  Let~$\cC_V$ be the subalgebra in~$\Phi_N$
generated by the elements of the $n$-dimensional plane $P\subset F^*$.
In other words, the algebra $\cC_V$ is the image of the induced mapping
$$
\Sym(E^*)\longrightarrow \Phi_N=\Sym(F^*)/\<\phi_i^2,\,i=1,\dots,N\>.
$$
The algebra $\cC_V$ has an obvious grading~$\cC_V=\cC^0_V\oplus
\cC^1_V\oplus \cC^2_V\oplus \cdots$ by degree of elements.

Suppose~$E=\h$ and $V$ is the collections of coroots~$h_\alpha$,
$\alpha\in\Delta_+$.  For $\lambda\in \h^*=E^*$, 
$p^*(\lambda)=\Theta(L_\lambda)\in F^*$ is the curvature 
form (see Proposition~\ref{pr:curvature}).
Then $\cC_V=\CC(X)$ is the algebra generated by the curvature forms
$\Theta(L_\lambda)$.

In general, we have the following result.

\begin{theorem}
\label{th:AV}
The dimension of the algebra $\cC_V$ is equal to the number
$\ind(V)$ of independent subsets in~$V$.  Moreover, the dimension of
the $k$-th component $\cC^k_V$ is equal to the number of
independent subsets~$S\subset V$ such that $k=N-|S|-\act(S)$.
\end{theorem}

We can also describe the algebra $\cC_V$  as a quotient of a
polynomial ring.  Let us say that a hyperplane $H$ in $E$ is a 
{\it $V$-essential hyperplane} if the elements of the subset 
$\{v_i,\, i=1,\dots,N \mid v_i \in H \}$
span the hyperplane $H$.  Obviously, an essential hyperplane is uniquely
determined by the subset of indices $I_H=\{i\in\{1,\dots,N\}\mid v_i\not\in
H\}$.  We will call such subset $I_H$ a {\it $V$-essential index subset}.
Denote by $d(H)=d_V(H)=|I_H|$ the number of its elements.  
A nonzero vector $\lambda\in
E^*$ determines the hyperplane $H=\{x\in E\mid \lambda(x)=0\}$ in $E$.  Vectors
$\lambda_H\in E^*$ corresponding to essential hyperplanes $H$ will be called
{\it $V$-essential vectors}.  They are defined up to a nonzero factor.

\begin{theorem}
The algebra $\cC_V$ is naturally isomorphic to the quotient of the 
polynomial ring $\mathrm{Sym}(E^*) / \I_V$, 
where the ideal $\I_V$ is generated by the powers $(\lambda_H)^{d(H)+1}$
of $V$-essential vectors for all $V$-essential hyperplanes $H$ in $E$.
The isomorphism is induced by the embedding $p^*:E^*\to F^*$.
\label{th:Rel} 
\end{theorem}

\begin{remark} \ {\rm 
There are several equivalent definitions of essential subsets, as follows:

\noindent
{\bf 1.} \ 
An index subset $I=\{i_1,\dots,i_k\}\subset\{1,\dots,N\}$ is {\it
$V$-essential} if and only if the following two conditions are satisfied: {\rm
(i)} the coordinate plane $\left<\phi_{i_1},\dots,\phi_{i_k}\right>$ in $F^*$
has one-dimensional intersection with the plane $P$; {\rm({i}{i})} there is no
proper subset in $I$ that satisfies the condition {\rm (i)}.
For an $V$-essential hyperplane $H$, the vector $p^*(\lambda_H)\in F^*$ spans
the one-dimensional intersection of $P$ and the coordinate plane associated
with $I_H$.

\noindent
{\bf 2.} \
Let $\theta_1,\dots,\theta_n$ be any basis in $P$. 
The $V$-essential subsets are in one-to-one correspondence with the 
cycles in the vector set
$\{\phi_1,\dots,\phi_N,\theta_1,\dots,\theta_n\}$.  For a cycle
$\{\phi_{i_1},\dots,\phi_{i_k}, \theta_{j_1},\dots,\theta_{j_k}\}$, the subset
$\{i_1,\dots,i_k\}$ is $V$-essential.  Moreover,
every $V$-essential subset is of this form.  
} 
\label{rm:equiv}
\end{remark}

Note that the decomposition of the Grassmannian $G(n,N)$ of all $n$-dimensional
planes $P\subset F^*$  into strata with the same collection of essential subsets
coincides with the decomposition of $G(n,N)$ into small cells of
Gelfand-Serganova~\cite{GelS} because any two $P_1,\, P_2\in G(n,N)$ 
with the same collection of essential subsets have the same dimensions
of intersections with all coordinate subspaces.
Equivalently, $P_1$ and $P_2$
are in the same strata if and only if the corresponding collections of vectors
$V_1$ and $V_2$ define the same matroid, i.e., have the same collection of
independent subsets.

Let us apply the Theorem~\ref{th:Rel} to $\CC(X)$.
Let $\omega_1,\omega_2,\dots,\omega_l$ be the fundamental weights.  They
generate the weight lattice~$\hat{T}$.  Also let $d_i$ be the number of
positive roots $\alpha\in\Delta_+$ such that $\alpha(\omega_i)\ne 0$.

\begin{corollary}  
The algebra $\CC(X)$ is naturally isomorphic to the quotient
of the polynomial ring $\Sym(\h^*)/\mathcal{J}$,
the ideal~$\mathcal{J}\in\Sym(\h^*)$ is generated by the elements
$(w\cdot \omega_i)^{d_i+1}$, where $i=1,\dots,l$ and $w$ is an element of the
Weyl group~$W$.
This isomorphism is induced by the projection
$\Sym(\h^*)\to \CC(X)$ that send $\lambda\in\hat{T}$
to the curvature form $\Theta(L_\lambda)$. 
\end{corollary}

For type~$A$ flag manifolds this statement was earlier proved in~\cite{SS}.

\proof
The generators of the ideal $\mathcal{J}$, as described in
Theorem~\ref{th:Rel}, correspond to root subsystems in $\Delta$ of
codimension~$1$. Every subsystem of codimension~$1$ is congruent,
modulo the Weyl group, to the root subsystem in the hyperplane $H_i$ of
zeroes of a fundamental weight~$\omega_i$. Moreover, $d(H_i)=d_i$ is exactly
the number of positive roots that do not belong to this hyperplane. 
\endproof

\section{Dual Statement and Proof of Theorems~\ref{th:main} and~\ref{th:AV}}

Next we give a dual version of Theorem~\ref{th:AV}.  Let~$\Sym(E)$ be the
symmetric algebra of the vector space~$E$.  For a subset $S\subset V$, denote
$m(S) =\prod_{v\in S} v\in\Sym(E)$, a square-free monomial in~$v_i$'s.  Let
us denote by $\S_V$ a subspace in $\Sym(E)$ spanned by $2^N$ square-free
monomials $m(S)$, where $S$ ranges over all subsets in $V$. 
Let $\S^k_V$ be the $k$-th graded component of~$\S_V$.

\begin{theorem}
\label{th:independent}
The dimension of the subspace $\S_V$ in $\Sym(E)$ is equal to the number
$\ind(V)$ of independent subsets in~$V$.
Moreover, the dimension of~$\S^k_V$ is equal to the number of independent 
subsets $S\subset V$ such that $k = N-|S|-\act(S)$.
\end{theorem}

Let us choose a basis $x^1,\dots,x^n$ in $E$.  Let $x_1,\dots,x_n$ be
the dual basis in $E^*$, and let $\theta_i=p^*(x_i)$, $i=1,\dots,n$, be the
corresponding basis in $P$.

For $I=(1\leq i_1<i_2<\dots<i_k\leq N)$ and $J=(1\leq j_1\leq j_2 \leq  \cdots
\leq j_k\leq n)$, we denote by $a_{IJ}$ the coefficient of the square-free
monomial $\phi_{i_1}\cdots \phi_{i_k}\in \Phi_N$ in the product of generators
$\theta_{j_1}\cdots \theta_{j_k}\in \cC_V$.  Let $A_k=(a_{IJ})$ be the
${ N \choose k} \times {n+k-1\choose k}$-matrix formed by the~$a_{IJ}$.  
Clearly, $\dim
\cC^k_V$ is equal to the rank of the matrix $A_k$.

On the other hand, $a_{IJ}$ is also the coefficient of the monomial
$x^{j_1}\cdots
x^{j_k}\in \Sym(E)$ in the square-free monomial $v_{i_1}\cdots v_{i_k}\in
\S(V)$.  Thus $\dim \S^k_V$ is equal to the rank of the transposed matrix
$A_k^T$, which is the same as the rank of~$A_k$.  
We proved the following statement, which implies the
equivalence of Theorems~\ref{th:AV} and~\ref{th:independent}.

\begin{lemma} 
For $k=1,\dots,N$, we have $\dim\cC^k_V=\dim\S^k_V$.
\end{lemma}

Let us say that a subset $S\subset V$
is {\it robust\/} if there is no cycle~$C\subset V$ with minimal element~$v$
such that $S\cap C=\{v\}$.
\begin{lemma}
\label{le:NBC}
The number of robust subsets in~$V$ is equal to the number of independent
subsets.
Moreover, the number of $k$-element robust subsets is equal to
the number of independent subsets $S$ with $k=N-|S|-\act(S)$.
\end{lemma}

\proof
We present an explicit bijection between robust subsets and independent
subsets.  A subset~$A$ is the complement to a robust subset in~$V$
if and only if for any cycle $C$ with minimal element $v$, inclusion
$C\setminus\{v\}\subset A$ implies $C\subset A$.  We will call such subsets
{\it antirobust.}

For an independent subset $S$, let $M$ be the collection of all
externally active $v\in V\setminus S$.  Then~$S\cup M$ is antirobust.
Conversely, for an antirobust subset $A$, let $M$ be the collection of all
$v\in V$ such that $v$ is a minimal element in some cycle $C\subset A$.  Then
$A\setminus M$ is an independent subset.

Clearly, both these mapping are inverse to each other and the statement of
lemma follows.
\endproof

\begin{theorem}
\label{th:robust}
The set of square-free monomials~$m(S)$, where $S$ ranges over
robust subsets, forms a basis of the subspace $\S_V$.
\end{theorem}

First, we prove a weaker version of Theorem~\ref{th:robust}.

\begin{lemma}
\label{le:robust}
The square-free monomials~$m(S)$, where $S$ ranges over
robust subsets, span $\S_V$.
\end{lemma}

\proof Suppose not.  Let $m(R)$ be the maximal in the lexicographical order
square-free monomial which cannot be expressed linearly via the
monomials~$m(S)$ with robust~$S$.  Then there is a cycle
$C=\{v_{i_1},v_{i_2},\dots,v_{i_l}\}$ with the minimal element~$v=v_{i_1}$ such
that $R\cap C=\{v\}$.  We can replace $v$ in the monomial $m(R)$ by a linear
combination of $v_{i_2}, v_{i_3},\dots,v_{i_l}$.  Thus $m(R)$ is a linear
combination of square-free monomials which are greater than $m(R)$ in the
lexicographical order.  By assumption each of these monomials can be expressed
via the monomials $m(S)$ with robust $S$.  Contradiction.  \endproof

We can now conclude the proof.

\proofof{Theorems~\ref{th:independent} and~\ref{th:robust}}
Lemmas~\ref{le:NBC} and~\ref{le:robust} imply the inequality $\dim
\S_V\leq \ind(V)$.
In view of these two lemmas it is enough to show that $\dim \S_V$ is
actually equal
to $\ind(V)$.

We prove this statement by induction on $|V|$.  If the linear span of
vectors in~$V$
is one-dimensional, then both $\dim \S_V$ and $\inv(V)$ are equal to
the number of non-zero vectors in~$V$ plus one.  This establishes the base
of induction.

Assume that $v=v_N$ is a nonzero vector in~$V$.
Let $V'=V\setminus\{v\}=\{v_1,\dots,v_{N-1}\}$, and
let $V''$ be the collection of images of vectors $v_1, \dots, v_{N-1}$
in the quotient space $E/\<v\>$.  It follows from the definition of
independent subset that $\ind(V) = \ind(V') + \ind(V'')$.
Assume by induction that $\dim \S_{V'}=\ind(V')$ and
$\dim \S_{V''}=\ind(V'')$.

Clearly, $\S_V$ is spanned by $\S_{V'}+\S_{V'}v$.   Both $\S_{V'}$ and 
$\S_{V'}v$
have same dimensions.  Hence, $\dim \S_V= 2 \dim \S_{V'}-\dim (\S_{V'}\cap
\S_{V'}v)$.  Let $\pi:\S_{V'}\to \S_{V''}$ be the natural projection.  Then
$\S_{V'}\cap \S_{V'}v\subset \Ker(\pi)$.  Thus $\dim \S_{V'} - \dim (\S_{V'}\cap
\S_{V'}v) \geq  \dim \S_{V''}$ and $\dim \S_{V}\geq \dim \S_{V'}+\dim \S_{V''}
=\ind(V')+\ind(V'')=\ind(V)$.  Coupled with the inequality $\dim \S_{V}\leq
\ind(V)$, this produces the required statement.

This finishes the proof of Theorems~\ref{th:independent} and~\ref{th:robust}
and thus of Theorems~\ref{th:main} and~\ref{th:AV}.
\endproof

\section{Proof of Theorem~\ref{th:Rel}}

Let $\Ess_V$ denote the set of all $V$-essential  hyperplanes in $E$.
Recall that $\I_V$ is the ideal in $\Sym(E^*)$ generated by the powers
of essential vectors
$(\lambda_H)^{d(H)+1}$, $H\in \Ess_V$ (see
Theorem~\ref{th:Rel}).  The embedding $p^*:E^*\to F^*$ induces the mapping
$\Sym(E^*)\to\Phi_N$, whose image is the algebra $\cC_V$.  Let $\tI_V$ denote
the kernel of this mapping, which we will call the {\it vanishing ideal} of
$V$.  Theorem~\ref{th:Rel} amounts to the identity of ideals $\I_V=\tI_V$.

The proof relies on a couple of simple lemmas.  As in the previous section we
assume that $v_N$ is a nonzero vector in~$V$.  Let
$V'=V\setminus\{v_N\}=\{v_1,\dots,v_{N-1}\}$, and let $V''$ be the collection
of images of vectors $v_1, \dots, v_{N-1}$ in the quotient space $E/\<v_N\>$.
We also denote by $\Ess_{V'}$ and $\Ess_{V''}$ the sets of $V'$-essential and
$V''$-essential hyperplanes in the corresponding spaces.

The dimension of the span of vectors in $V''$ is $n-1$.  The dimension of the
span of $V'$ can be either $n-1$ or $n$.  For a hyperplane $H$ in  $E/\<v_N\>$,
let $\overline{H}=H\oplus \<v_N\>$ be the hyperplane in $E$.  Also for a
collection of hyperplanes $C$ in $E/\<v_N\>$, let $\overline{C}= \{\overline{H} 
\mid H \in C\}$ be the collection of hyperplanes in $E$.

\begin{lemma}
\label{le:IND} 
{\rm (a)}  If $\dim V' = n-1$ then 
$\Ess_{V} = \{\<V'\>\} \cup \overline{\Ess_{V''}}$. 

\noindent
{\rm (b)} If $\dim V' =n$ then
$\Ess_V= \Ess_{V'} \,\bigcup\, \overline{\Ess_{V''}}$.
For $H\in (\Ess_{V'}\setminus\overline{\Ess_{V''}})$,
we have $d_{V}(H)=d_{V'}(H)+1$.  For $H\in \Ess_{V''}$,
we have $d_V(\overline{H}) = d_{V''}(H)$.
\end{lemma}

\proof 
The part (a) is left as an easy exercise for the reader.
In order to prove (b) we first assume that a hyperplane $H$ 
contains $v_N$.  Then $H$ is $V$-essential if and only if its 
projection $H/\<v_N\>$ is $V''$-essential.
Suppose that a hyperplane $H$ does not contain $v_N$.
Then $H$ is $V$-essential if and only if it is $V'$-essential and
its projection $H/\<v_N\>$ is not $V''$-essential.
The equalities for the numbers $d(H)$ are also obvious.
This proves the lemma.
\endproof

Recall that the collection of vectors $V$ is associated with a plane $P\in
G(n,N)$.
Let $P'$ and $P''$ be the planes associated with vector sets $V'$ and $V''$,
respectively.  Namely, $P'$ is the projection of $P$  along $\phi_N$
onto the hyperplane
$\mathcal{H}=\<\phi_1,\dots,\phi_{N-1}\>$ 
spanned by the first $N-1$ coordinate vectors;
and $P''$ is the intersection of $P$ with the same hyperplane~$\mathcal{H}$.

We can choose the basis $x_1,\dots,x_n$ in $E^*$ 
and the corresponding basis $\theta_i=p^*(x_i)$, $i=1,\dots,n$ in $P$ such that
$\theta_1,\dots,\theta_{n-1}\in \mathcal{H}$ and 
$\theta_n\in \phi_N+\mathcal{H}$.
Also denote $\ttheta_{n}=\theta_n-\phi_N \in \mathcal{H}$.
Then $\theta_{1},\dots,\theta_{n-1},\ttheta_n$ span the projected space 
$P'$.  The space $\Sym(E^*)$ can be identified with the polynomial ring
$\C[x_1,\dots,x_n]$

\begin{lemma}
\label{le:DER} The vanishing ideal $\tI_V$
consists of all polynomials $f\in\C[x_{1},\dots,x_{n}]$
such that both $f$ and $\pa f/\pa x_n$ belong to
the vanishing ideal $\tI_{V'}$. In
particular, $\tI_{V}\subseteq\tI_{V'}.$
\end{lemma}

\proof Recall that the ideal $\tI_{V}$ consists of
all polynomials in $x_{1},\ldots,x_{n}$ which vanish in the algebra 
$\Phi_N$ upon substituting of the $\theta_i$ instead of the $x_{i}$.
Taylor's expansion in the algebra $\Phi_N$ gives
$$
f(\theta_{1},\dots,\theta_{n})=f(\theta_{1},\dots,\ttheta_{n}+\phi_{N})=
f(\theta_{1},\dots,\ttheta_{n})+
({\pa f / \pa x_n}) (\theta_{1},\dots,\theta_{n-1},\ttheta_{n}) \phi_{N}.
$$
We have only two non-vanishing terms in the right hand side, since $\phi_N^2=0$.

The polynomial $f(x_{1},\dots,x_{n})$ belongs to $\tI_{V}$ if and only if
both $f(\theta_{1},\dots,\ttheta_{n})$ and
$({\pa f / \pa x_n}) (\theta_{1},\dots,\theta_{n-1},\ttheta_{n})$ vanish 
in the algebra $\Phi_N$.
\endproof


We now conclude the proof of the identity $\I_V=\tI_V$.  The inclusion
$\I_V\subseteq\tI_V$ is straightforward.  Indeed, every $p^*(\lambda_H)$ is a
linear combination of $d(H)$ different $\phi_i$'s.  Thus $(\lambda_H)^{d(H)+1}$
maps to zero in the algebra~$\Phi_N$.

We prove the identity $\I_V=\tI_V$ by induction on $|V|$.  
The base of induction is the trivial case $V=\{v_1\}$. 
For $|V|\geq 2$, assume by induction that the statement is true for both
$V'$ and $V''$.  Take any $f\in\tI_{V}$.  Our goal is to show that $f\in \I_V$.

By Lemma~\ref{le:DER}, the polynomial ${\pa f}/{\pa x_{n}}$ belongs 
to $\tI_{V'}$.  
Notice that for any $H\in \Ess_{V'} \cap \overline{\Ess_{V''}}$ 
the coordinate expansion of the corresponding 
$V'$-essential vector $\lambda_{H}$ does not involve $x_{n}$.  
Indeed, any $H\in \overline{\Ess_{V''}}$ contains $v_N$, thus
$\lambda_H(v_N)=0$.  On the other hand,
$x_1(v_N)=\cdots=x_{n-1}(v_N)=0$, and $x_n(v_N)=1$.

By inductive assumption, one has
$$
{\pa f}/{\pa x_{n}}=\sum_{H\in \Ess_{V'} \setminus \overline{\Ess_{V''}}} 
p_{H}\, \lambda_H^{d_{V'}(H)+1}+
\sum_{H\in \Ess_{V'} \cap \overline{\Ess_{V''}}} 
p_{H}\, \lambda_H^{d_{V'}(H)+1},
$$
where the $p_{H}$ are certain polynomials in $x_{1},\dots,x_{n}$.
The $\lambda_H$ in the second sum do not involve $x_n$.
Thus, integrating the above expression with respect to $x_n$, one 
deduces that there exists a polynomial $\bar f$ of the form
$$
\bar f=\sum_{H\in \Ess_{V'} \setminus \overline{\Ess_{V''}}} 
\bar p_{H}\, \lambda_H^{d_{V'}(H)+2}+
\sum_{H\in \Ess_{V'} \cap \overline{\Ess_{V''}}} \bar p_{H} \,
\lambda_H^{d_{V'}(H)+1}
$$
satisfying ${\pa \bar f}/{\pa x_{n}}={\pa f}/{\pa x_{n}}.$
By Lemma~\ref{le:IND}, the polynomial $\bar f$ can be written as
$$
\bar f=\sum_{H\in\Ess_{V}} \bar p_H\, \lambda_H^{d_V(H)+1}
$$
and thus belongs to  $\I_{P}$. The difference
$\hat f=f-\bar f$ belongs to $\tI_{P}$ and is independent
on $x_{n}$.  Thus $\hat f\in {\tI}_{V''}$. 
By induction hypothesis, $\hat f\in \I_{V''}\subset \I_V$.
Thus $f=\bar f + \hat f\in\I_V$.
The statement follows. \hfill Q.E.D.

\section{Remarks and Open Problems}

The algebra $\CC(X)$ for type~A flag manifolds $X=SL(n)/B$ was studied in
more details in~\cite{SS} and~\cite{PSS} motivated by~\cite{Ar}.  In this case, we first proved
Theorem~\ref{th:main} in~\cite{PSS} using a different approach based on a
presentation of~$\CC(X)$ as a quotient of a polynomial ring 
(cf.~Theorem~\ref{th:Rel}).
The theorem claims that the dimension of~$\CC(X)$ is equal to the
number forests on~$n$ labelled vertices whereas the dimension of~$\CC^k(X)$ is
equal to the number of forests with ${n \choose 2} - k$ inversions.  The
statement about forests was initially conjectured in~\cite{SS} and the
statement concerning inversions was then guessed by R.~Stanley.  In~\cite{PSS}
we also discuss various generalizations of the ring~$\CC(X)$.

A natural open problem is to extend the results to homogeneous manifolds
$G/P$, where~$P$ is a parabolic subgroup.
Formulas for curvature forms on~$G/P$ can be found in \cite{GS}.

It is also intriguing to find the links between the algebra~$\CC(X)$ and
the arithmetic Schubert calculus, see H.~Tamvakis~\cite{Ta1, Ta2}.

B.~Kostant pointed out that the algebra~$\CC(X)$ is related to the
$S(\g)$-module $\wedge\g$ studied in~\cite{K}.  It would be interesting to
investigate this relationship.

%
%
%
%

\bigskip
The authors are  grateful to Vladimir Igorevich Arnold, Richard Stanley, 
Alek Vainshtein, and Andrei Zelevinsky for stimulating discussions and 
helpful suggestions.

\bigskip

\end{document}